\input amstex
\documentstyle{amsppt}
\topmatter
\NoRunningHeads
\title
On the spectral asymptotics in domains with long boundary
\endtitle
\author
Leonid Friedlander
\endauthor
\affil
University of Arizona
\endaffil
\email
friedlan@math.arizona.edu
\endemail
\abstract
I discuss a simple toy problem for the Dirichlet Laplacian  in a sequence of domains where the contribution of the boundary to the spectral asymptotics
is of the same order as the contribution from the interior.

\endabstract
\endtopmatter
\document
\head 1. Introduction
\endhead

This note is motivated by some remarks in [DMO]. Dupius, Mazo, and Onsager studied the surface contribution to the specific heat in elasticity; that amounts to going to the
second term in the Weyl eigenvalue asymptotics for the Lame operator in a bounded domain. They discussed possible physical scenarios when the contribution of the surfase boundary
is noticeable, i.e. comparable to the contribution of the interior. This may happen when the boundary is rough, and its area is big enough.

Let me describe a toy two-dimensional problem that I will investigate. We start with the square $Q=[0,1]\times [0,1]\subset \Bbb R^2$. Let $q$ be a positive, integer number, and let
$T_q=[0, 1/(2q)]\times [0,h]$ be a rectangle; $h$ is a fixed number. Denote by $T_{q,j}=T_q+((j-1)/q, 1)$, $j=1,\ldots, q$, the parallel shifts of $T_q$ by the vectors $((j-1)/q, 1)$.
I will study the eigenvalue asymptotics for the Dirrichlet Laplacian in
$$D_q=Q\cup\bigl(\cup_{j=1}^q T_{q,j}\bigr).
$$
The area of $D_q$ and the circumference of $D_q$ equal $1+(h/2)$ and $2hq+4$, respectively. We will be interested in large values of $q$.

Let us recall the two term Weyl spectral asymptotics for planar domains.
Let $\Omega$ be a bounded planar domain with a piecewise smooth boundary, and let $N_\Omega(\lambda)$ be the counting function of the Dirichlet spectrum of
$\Omega$; that is the number of eigenvalues of the Dirichlet Laplacian in $\Omega$, which are smaller than or equal to $\lambda$. Under a "small number of periodic billiard trajectories in $\Omega$" assumption (e.g., see [SV]), one has
$$N_\Omega(\Lambda)=\frac{A(\Omega)}{4\pi}\lambda-\frac{L(\Omega)}{4\pi}\sqrt{\lambda}+o(\sqrt{\lambda}),\quad \lambda\to\infty,\tag 1$$
where $A(\Omega)$ is the area of $\Omega$ and $L(\Omega)$ is its circumference. The first and the second term in (1) are comparable if
$\lambda$ is of the order of $(L(\Omega)/A(\Omega))^2$; if applied to $D_q$ that means $\lambda\sim q^2$. Therefore, I set
$$\lambda=\mu q^2,\tag 2$$
and I will look at what happens when $q\to\infty$. Here $\mu$ is a constant; we will keep track of the $\mu$-dependence. Formally substituting $A(D_q)$ and $L(D_q)$ in (2) one gets
$$ N_{q,weyl}(\lambda)\sim \biggl(\frac{2+h}{8\pi}\mu-\frac{h}{2\pi}\sqrt{\mu}\biggr)q^2=c_{weyl}(\mu)q^2.\tag 3$$
As we will see, the actual asymptotics of $N_{D_q}(\mu q^2)$ is somewhat different. The main result of this note is
\proclaim{Theorem}
One has
$$N_{D_q}(\mu q^2)=c(\mu)q^2+O(q)$$
where
$$c(\mu)=\frac{\mu}{4\pi}+\frac{h}{\pi}\sqrt{\mu}\sum_{l=1}^{[\sqrt{\mu}/2\pi]}\sqrt{1-\frac{4\pi^2l^2}{\mu}}.
$$
Here $[\cdot]$ is the floor of a number.
\endproclaim
\head 
2. The Dirichlet-to-Neumann operator
\endhead
The Dirichlet-to-Neumann operator is the tool that I will use for estimating $N_{D_q}(\lambda)$. Let me start from a bounded domain $\Omega$ in $\Bbb R^d$ with a piecewise smooth boundary.
 Let $\Gamma=\partial\Omega$ be the boundary of $\Omega$, and let $\Sigma\subset\Gamma$ be a domain in $\Gamma$. In addition to the Dirichlet problem for the Laplacian
in $\Omega$, we will also consider the mixed spectral problem (sometimes it is called the Zaremba problem):
$$\Delta u+\lambda  u=0\ \text{in $\Omega$},\quad u=0\ \text{in $\Gamma\setminus\Sigma$},\quad \frac{\partial u}{\partial\nu}=0\ \text{in $\Sigma$}.$$
Here $\partial/\partial\nu$ is the outward normal derivative. In the case $\Sigma=\Gamma$, this is the Neumann problem. It has discrete spectrom: its eigenvalues are bracketted between eigenvalues of the Neumann problem
and the eigenvalues of the Dirichlet problem for the Laplacian. Let $N_{\Omega,\Sigma}(\lambda)$ be the counting function of the spectrum of this problem; $N_{\Omega}(\lambda)$ is,
as above, the counting function of the spectrum of the Dirichlet problem.

Now I introduce the Dirichlet-to-Neumann operator $R_\Sigma(\lambda)$. It is defined for values of $\lambda$ that are outside of the Dirichlet spectrum of the Laplacian, it acts on functions  on $\Sigma$, and its action is given by the formula
$$R_\Sigma(\lambda)f=\frac{P_\Sigma(\lambda)f}{\partial\nu}\bigg\vert_\Sigma \tag 4 $$
where $u(x)=(P_\Sigma f)(x)$ is the solution to the problem
$$\Delta u+\lambda u=0\ \text{in $\Omega$},\quad u=0\ \text{in $\Gamma\setminus\Sigma$},\quad u(x)=f(x)\ \text{in $\Sigma$}.\tag 5$$
The operator $P_\Sigma(\lambda)$ is called the Poisson operator. In the case $\Sigma=\Gamma$, one gets the usual Dirichlet-to-Neumann operator $R(\lambda)$. If $\Gamma$ is smooth then $R(\lambda)$ is an elliptic pseudodifferential operator of order $1$. We will always assume $\lambda$ to be real; then the operator $R(\lambda)$ is self-adjoint and it is bounded from below.
The same is true for the operator $R_\Sigma(\lambda)$. One has to take a little bit more care defining $R_\Sigma(\lambda)$. The easiest way is to do it via the quadratic form
$$Q_\Sigma(f)=\int_\Omega\bigl(|\nabla u(x)|^2-\lambda |u(x)|^2\bigr)dx\quad\text{where}\quad u(x)=(P_\Sigma f)(x).\tag 6
$$
Let $n_\Sigma(\lambda)$ be the number of non-positive eigenvalues of $R_\Sigma(\lambda)$. The first basic formula that I will use is
$$N_{\Omega,\Sigma}(\lambda)-N_\Omega(\lambda)=n_\Sigma(\lambda).\tag 7$$
It was proved in [F] in the case $\Sigma=\Gamma$; the proof extends almost verbatim to an arbitrary $\Sigma$.

The second formula that I will use is an extension of (7), and its proof is also the same as the proof in [F]. Let $\Omega$ be divided by a hypersurface $\Sigma$ into two
disjoint parts, $\Omega_\pm$. The hypersurface $\Sigma$ is  a part of the boundary for both parts of $\Omega$. Suppose that $\lambda$ is a real number, outside of the spectrum of the Dirichlet Laplacian in both $\Omega_\pm$. Let $R^\pm_\Sigma(\lambda)$ be the Dirichlet-to-Neumann operators related to $\Omega_\pm$, and $R_\Sigma(\lambda)=R^+_\Sigma(\lambda)+R^-_\Sigma(\lambda)$. Notice that outward directions to $\Omega_\pm$ are opposite to each other. If, as above, $n_\Sigma(\lambda)$ is the number of non-positive eigenvalues of $R_\Sigma(\lambda)$  then
$$
N_\Omega(\lambda)=N_{\Omega^+}(\lambda)+N_{\Omega^-}(\lambda)+n_\Sigma(\lambda).\tag 8
$$
\head
3. Proof of the theorem and remarks
\endhead
Let 
$$
\Sigma_q=\bigl[\cup_{j=0}^{q-1} [j/q, (2j+1)/(2q)]\bigr]\times \{1\}.
$$
Then $\Sigma_q$ divides $D_q$ into $D_q^+=Q$ and $D_q^-=\cup _{j=1}^q T_{q,j}$. Let $R_q=R_q^++R_q^-$ be the Dirichlet-to-Neumann operator that corresponds to $\Sigma_q$. Then
$$N_{D_q}(\mu q^2)=N_{Q}(\mu q^2)+N_{\cup T_{j,q}}(\mu q^2)+n_{\Sigma_q}(\mu q^2)\tag 9$$
(see (8)).  Firstly, I will show that
$$n_{\Sigma_q}(\mu q^2)=O(q).\tag 10$$
Indeed, the number of non-positive eigenvalues for a sum of two operators does not exceed the sum of numbers of non-positive eigenvalues of both of them. Therefore,
$$
n_{\Sigma_q}(\mu q^2)\leq n_{Q,\Sigma_q}(\mu q^2)+q n_{T_q,B_q}(\mu q^2)
$$
where $B_q-[0,1/(2n)]\times \{0\}$ is the base of the rectangle $T_q$ and $n_{\cdot}$ is the number of non-positive eigenvalues of the operator $R_{\cdot}$.
The equation (7) implies
$$ n_{Q,\Sigma_q}(\mu q^2)=N_{Q,\Sigma_n}(\mu q^2)-N_Q(\mu q^2).$$
Notice that $N_{Q,\Sigma_q}(\mu q^2)$ does not exceed the counting function of eigenvalues of the Neumann Laplacian, $\Cal N_Q(\mu q^2)$. Then,
$\Cal N_Q(\mu q^2)-N_Q(\mu q^2)=O(q)$; so
$$ n_{Q,\Sigma_q}(\mu q^2)=O(q).$$
The eigenfunctions of the operator $R_{T_q,B_q}$ are $\sin (2qk\pi x)$; the corresponding eigenvalues are $-v_k'(0)/v_k(0)$ where $v_k(y)$ solves the equation
$$v_k''+(\mu q^2-4\pi^2q^2k^2)v=0,$$
with $v_k(h)=0$. When $k>\sqrt{\mu}/(2\pi)$, this eigenvalue is positive. Therefore, the number of non-positive eigevalues does not exceed $\sqrt{\mu}/(2\pi)$. This finishes the proof of (10).

One has
$$N_Q(\mu q^2)=\frac{\mu}{4\pi}q^2+O(q)\tag 11$$
and 
$$
N_{\cup T_{j,q}}(\mu q^2)=qN_{T_q}(\mu q^2)\tag 12.
$$
The eigenvalues of the Dirichlet Laplacian in $T_q$ are $4\pi^2 l^2q^2+(\pi^2/h^2)k^2$ where $l$ and $k$ are positive, integer numbers. Such an eigenvalue does not exceed $\mu q^2$ if 
$k\leq (qh\pi)\sqrt{\mu-4\pi^2l^2}$. Therefore,
$$
N_{T_q}(\mu q^2)=\sum_{l=1}^{[\sqrt{\mu}/(2\pi)]}\biggl[\frac{qh}{\pi}\sqrt{\mu-4\pi^2l^2]}\biggr]=\frac{qh}{\pi}\sum_{l=1}^{[\sqrt{\mu}/(2\pi)]}\sqrt{\mu-4\pi^2l^2}+O(1).\tag 13
$$
Combining formulas  (9)--(13) one gets the statement of the theorem.
\qed

Let us compare the constants $c(\mu)$ in the theorem and $c_{weyl}(\mu)$ in (3). Let
$$f(x)=\sqrt{1-\frac{4\pi^2x^2}{\mu}}\ \ \text{and}\ \ m=\biggl[\frac{\sqrt{\mu}}{2\pi}\biggr].$$
If $\mu\geq 4\pi^2$ then one can apply the Euler--Maclaurin formula,
$$\sum_{l=1}^m f(x)dx=\int_0^m f(x)dx+\frac{f(m)}{2}-\frac{1}{2}+\int_0^m \biggl(\{x\}-\frac{1}{2}\biggr)f'(x)dx\tag 14
$$
where $\{x\}=x-[x]$ is the fractional part of $x$. Notice that $f(0)=1$. Then
$$\int_0^{\sqrt{\mu}/(2\pi)}f(x)dx=\frac{\pi\sqrt{\mu}}{8\pi},$$
so
$$c(\mu)=c_{weyl}(\mu)+\frac{h}{\pi}\sqrt{\mu}\biggl( \frac{f(m)}{2}-\int_m^{\sqrt{\mu}{2\pi}}f(x)dx+\int_0^m \biggl(\{x\}-\frac{1}{2}\biggr)f'(x)dx\biggr).\tag 15
$$
Notice that the function $f(x)$ is positive, decreasing, and concave on the interval $[0.\sqrt{\mu}/2(\pi)]$: that makes both quantities
$$\frac{f(m)}{2}-\int_m^{\sqrt{\mu}{2\pi}}f(x)dx\quad\text{and}\quad \int_0^m \biggl(\{x\}-\frac{1}{2}\biggr)f'(x)dx$$
negative. We conclude that
$$
c(\mu)<c_{weyl}(\mu).\tag 16
$$
If $\mu<4\pi^2$ then the sum in the expression $c(\mu)$ drops, and only the interior of the square cntributes to the asymptotic. The relation (16) still holds for $\mu>16$, and, for
$\mu<16$, the opposite inequality $c(\mu)>c_{weyl}(\mu)$ takes place.

\enddocument